\documentclass{article}

\usepackage[dvips]{graphicx}
\usepackage{amssymb}
\usepackage{amsmath}

\usepackage{caption}
\usepackage{hyperref}
\usepackage{arcs}
\usepackage {xcolor} 
\usepackage{subfigure}
\usepackage{authblk}
\usepackage{calrsfs}

\newtheorem{theorem}{Theorem}

\usepackage{graphicx} % Required for inserting images

\title{Unisolvence of random Kansa collocation by Thin-Plate Splines for the Poisson equation}
\author[1]{F. Dell'Accio}
\affil[1]{University of Calabria, Rende (CS), Italy}
\author[2]{A. Sommariva} 
\affil[2,3]{University of Padova, Italy}
\author[3]{M. Vianello}
\date{\today}

\begin{document}

\maketitle

\begin{abstract}

Existence of sufficient conditions for unisolvence of Kansa unsymmetric collocation for PDEs is still an open problem. 
In this paper we make a first step in this direction, proving that unsymmetric collocation matrices with  Thin-Plate Splines for the 2D Poisson equation are almost surely nonsingular, when the discretization points are chosen randomly on domains with analytic boundary.  
\end{abstract}

\section{Introduction}

Kansa unsymmetric collocation, originally proposed in the mid '80s \cite{K86}, has become over the years a popular meshless method for the discretization of boundary value problems for PDEs. Despite its wide and successful adoption for the numerical solution of a variety of physical and engineering problems (cf. e.g. \cite{CFC14} with the references therein), a sound theoretical foundation concerning unisolvence of the corresponding linear systems is still missing. Indeed, it was shown by Hon and Schaback \cite{HS01} that there exist point configurations that lead to singularity of the collocation matrices, though these are very special and ``rare''cases. For this reason greedy and other approaches have been developed to overcome the theoretical problem and ensure invertibility, cf. e.g. \cite{LOS06,SW06}. On the other hand, in the textbook \cite{F07} one can read : ``{\em Since the numerical experiments by Hon and Schaback show that Kansa's method cannot be well-posed for arbitrary center locations, it is now an open question to find sufficient conditions on the center locations that guarantee invertibility of the Kansa matrix}'', and the situation does not seem to have changed so far. 

In this paper we make a first step in this direction, proving that unsymmetric collocation matrices with  Thin-Plate Splines (without polynomial addition) for the 2D Poisson equation are almost surely nonsingular, when the discretization points are chosen randomly on domains with analytic boundary.  
Though TPS are not the most adopted option for Kansa collocation, they have been often used in the meshless literature, cf. e.g. \cite{CFC14,CDN06,ZDC00} with the references therein. 
One of their most
relevant features is that they are scale invariant, thus avoiding the delicate matter of the scaling choice with scale dependent RBF, which is still an active research topic, cf. e.g. \cite{CDR20,LS23}.
On the other hand, the fact that TPS without polynomial addition can guarantee unisolvence in the interpolation framework has been recently recognized experimentally in \cite{P22} and theoretically in \cite{BSV23,DASV23}.

As we shall see, one of the key aspects is that Thin-Plate Splines $\phi(\|P-A\|_2)$, which correspond to the radial functions 
\begin{equation} \label{TPS}
\phi(r)=r^{2\nu} \log(r)\;,\;\;\nu\in \mathbb{N}\;,
\end{equation}
are {\em real analytic functions} off their center $A$, due to analyticity of the univariate functions $\log(\cdot)$ and $\sqrt{\cdot}$ in $\mathbb{R}^+$.
Analiticity together with the presence of a singularity at the center will be the key ingredients of our unisolvence result by random collocation. 

\section{Unisolvence of random Kansa collocation}
Consider the Poisson equation with Dirichlet boundary conditions (cf. e.g. \cite{E98})

\begin{equation}\label{poisson_eq}
\left\{
\begin{array}{l}
\Delta u(P)=f(P)\;,\;P\in \Omega
\\
u(P)=g(P)\;,\;P\in \partial \Omega=\gamma([a,b])\;,
\end{array} 
\right .
\end{equation}
where we assume that $\Omega\subset \mathbb{R}^2$ is a domain with analytic boundary (a bounded connected open set whose boundary is an analytic curve), namely a curve $\gamma:[a,b]\to \mathbb{R}^2\;,\;\gamma(a)=\gamma(b)$, that is {\em analytic} and {\em regular} (i.e. $\gamma'(t)\neq (0,0)$ for every $t\in [a,b]$).

In Kansa collocation (see e.g. \cite{F07,HS01,K86,LOS06,SW06,W05}) one determines a function
\begin{equation} \label{u_N}
u_N(P)=\sum_{j=1}^n{c_j\,\phi_j(P)}+\sum_{k=1}^m{d_k\,\psi_k(P)}\;,\;\;N=n+m\;,
\end{equation}
where
\begin{equation} \label{phij}
\phi_j(P)=\phi(\|P-P_j\|_2)\;,\;\;\{P_1,\dots,P_n\}\subset \Omega\;,
\end{equation}
\begin{equation} \label{psik}
\psi_k(P)=\phi(\|P-Q_k\|_2)\;,\;\;\{Q_1,\dots,Q_m\}\subset 
\partial\Omega\;,
\end{equation}
such that
\begin{equation}\label{poisson_disc}
\left\{
\begin{array}{l}
\Delta u_N(P_i)=f(P_i)\;,\;i=1,\ldots,n
\\
u_N(Q_h)=g(Q_h)\;,\;h=1,\ldots,m\;.
\end{array} 
\right .
\end{equation}

The following facts will be used below. Defining $\phi_A(P)=\phi(\|P-A\|)$, we have $\phi_A(B)=\phi_B(A)$ and 
$\Delta\phi_A(B)=\Delta\phi_B(A)$. In fact, the Laplacian in polar coordinates 
centered at $A$ (cf. e.g. \cite[Ch.2]{E98}) is the radial function
\begin{equation} \label{lapl}
\Delta\phi_A={\frac{\partial^2 \phi}{\partial^2 r}}+{\frac{1}{r}} \frac{\partial \phi}{\partial r}=
4\nu r^{2(\nu-1)}(\nu\log(r)+1)\;.
\end{equation}
Moreover, $\phi_A(A)=0$ and 
$\Delta\phi_A(A)=0$ for $\nu\geq 2$, since $\Delta\phi\to 0$ as $r\to 0$. 

Kansa collocation can be rewritten in matrix form as
\begin{equation} \label{Kansa-system}
\left(\begin{array}{cc}
\Delta\Phi & \Delta\Psi\\
\\
\Phi & \Psi
\end{array} \right) 
\left(\begin{array}{c}
\mathbf{c}\\
\\
\mathbf{d}
\end{array} \right) 
=
\left(\begin{array}{c}
\mathbf{f}\\
\\
\mathbf{g}
\end{array} \right) 
\end{equation}
where the block matrix is 
$$
K_N=K_N(\{P_i\},\{Q_h\})=\left(\begin{array}{cc}
\Delta\Phi & \Delta\Psi\\
\\
\Phi & \Psi
\end{array} \right) 
$$
$$
=\left(\begin{array} {ccccccc}
0 & \cdots &  \cdots & \Delta\phi_n(P_1)
& \Delta\psi_1(P_1) & \cdots & \Delta\psi_m(P_1) \\
\\
\vdots & \ddots  &  & \vdots & \vdots  & \cdots & \vdots\\

\vdots &  & \ddots & \vdots & \vdots  & \cdots & \vdots\\
\\
\Delta\phi_1(P_n) & \cdots & \cdots & 0
& \Delta\psi_1(P_n) & \cdots & \Delta\psi_m(P_n) \\
\\
\phi_1(Q_1) & \cdots & \cdots & \phi_n(Q_1)
& 0 & \cdots & \psi_m(Q_1) \\
\\
\vdots & \cdots & \cdots & \vdots & \vdots  & \ddots & \vdots\\
\\
\phi_1(Q_m) & \cdots & \cdots & \phi_n(Q_m)
& \psi_1(Q_m) & \cdots & 0 \\
\\
\end{array} \right)
$$
 and $\textbf{f}=\{f(P_i)\}_{i=1,\ldots,n}$, ${\textbf{g}}=\{g(Q_h)\}_{h=1,\ldots,m}$.

 We can now state and prove our main result.

\begin{theorem}
Let $K_N$ be the TPS Kansa collocation matrix defined above, with $N=n+m\geq 2$, where $\{P_i\}$ is a sequence of independent uniformly distributed random points in $\Omega$, and $\{Q_h\}$ a sequence of independent uniformly distributed points on $\partial\Omega$. Namely, $\{Q_h\}=\{\gamma(t_h)\}$ with $\{t_h\}$ sequence of independent identically distributed random abscissas in $(a,b)$ with respect to the arclength density $\|\gamma'(t)\|_2/L$, $L=length(\gamma([a,b]))$. 

Then for every $N\geq 2$ the matrix $K_N$ is a.s. (almost surely) nonsingular.
\end{theorem}

\vskip0.5cm 
\noindent
{\bf Proof.}
The proof proceeds by complete induction on $N$. 
For the induction base,
we prove that $\mbox{det}(K_N)$ is a.s. nonzero for $N=2$, that is for $n=2$ and $m=0$, or $n=0$ and $m=2$, or $n=1$ and $m=1$. In the first case, $$\mbox{det}(K_2)=-\Delta\phi_2(P_1)\Delta\phi_1(P_2)=-(\Delta\phi_1(P_2))^2
$$
$$
=-16\nu^2\|P_2-P_1\|_2^{4\nu-4}
\left(\nu\log(\|P_2-P_1\|_2)+1\right)^2$$
which vanishes iff $P_2=P_1$ (an event with null probability) or $P_2$ falls on (the intersection with $\Omega$ of) the curve   
$\nu\log(\|P-P_1\|_2)+1=0$, that is on the circle
$$\|P-P_1\|_2^2=\exp(-2/\nu)\;.
$$
But this event has null probability, since any  algebraic curve is a null set in $\mathbb{R}^2$.

In the second case, 
$$\mbox{det}(K_2)=-\psi_2(Q_1)\psi_1(Q_2)=-\psi_1^2(Q_2)
=-\psi_1^2(\gamma(t_2))\;.
$$
Now, given $P_1=\gamma(t_1)$, the function $\lambda(t)=\psi_1^2(\gamma(t))$ is analytic in $(a,t_1)$ and in $(t_1,b)$. Then $\psi_1^2(\gamma(t_2))$ is zero iff $t_2=t_1$ (an event that has null probability), or $t_2$ falls on the zero set of $\lambda$ 
in $(a,t_1)$ or $(t_1,b)$. Again this event has null probability since the zero set of an univariate analytic function in an open interval is a null set (cf. \cite{KP02,M20}).

As for the third case, 
assume that $Q_1$ is chosen on the boundary (randomly or not) and that $P_1$ is chosen randomly in the interior. Since 
$$
\mbox{det}(K_2)=-\phi_1(Q_1)\Delta\psi_1(P_1)
$$
$$
=-4\nu\phi_1(Q_1)
\|P_1-Q_1\|_2^{2\nu-2}
\left(\nu\log(\|P_1-Q_1\|_2)+1\right)\;,
$$
and $\phi_1(Q_1)\neq 0$ being $P_1\neq Q_1$, the determinant vanishes if and only if $P_1$ falls on (the intersection with $\Omega$ of) the curve   
$\nu\log(\|P-Q_1\|_2)+1=0$, that is on the circle
$$\|P-Q_1\|_2^2=\exp(-2/\nu)
$$
and again this event has null probability.

For the inductive step, we consider separately the case where a boundary  point is added, for which we define the matrix 
$$
U(P)=\left(\begin{array} {cccccccc}
0 & \cdots &  \cdots & \Delta\phi_n(P_1)
& \Delta\psi_1(P_1) & \cdots & \Delta\psi_m(P_1) & \Delta\phi_1(P)\\
\\
\vdots & \ddots  &  & \vdots & \vdots  & \cdots & \vdots & \vdots\\
\vdots &  & \ddots & \vdots & \vdots  & \cdots & \vdots & \vdots\\
\\
\Delta\phi_1(P_n) & \cdots & \cdots & 0
& \Delta\psi_1(P_n) & \cdots & \Delta\psi_m(P_n) & \Delta\phi_n(P)\\
\\
\phi_1(Q_1) & \cdots & \cdots & \phi_n(Q_1)
& 0 & \cdots & \psi_m(Q_1) & \psi_1(P)\\
\\
\vdots & \cdots & \cdots & \vdots & \vdots  & \ddots & \vdots & \vdots\\
\\
\phi_1(Q_m) & \cdots & \cdots & \phi_n(Q_m)
& \psi_1(Q_m) & \cdots & 0 & \psi_m(P)\\
\\
\phi_1(P) & \cdots & \cdots & \phi_n(P)
& \psi_1(P) & \cdots & \psi_m(P) & 0\\
\\
\end{array} \right)
$$

Observe that in this case $K_{N+1}=U(Q_{m+1})$. Indeed, $\psi_k(Q_h)=\psi_h(Q_k)$ and $\Delta\phi_i(Q_{m+1})=\Delta\psi_{m+1}(P_i)$. 

Differently, if an interior point is added, we define the matrix

$$
V(P)=\left(\begin{array} {cccccccc}
0 & \cdots &  \cdots & \Delta\phi_n(P_1) & \Delta\phi_1(P)
& \Delta\psi_1(P_1) & \cdots & \Delta\psi_m(P_1) \\
\\
\vdots & \ddots  &  & \vdots & \vdots  & 
\vdots & \cdots & \vdots \\
\vdots &  & \ddots & \vdots & \vdots  & \vdots & \cdots & \vdots\\
\\
\Delta\phi_1(P_n) & \cdots & \cdots & 0
& \Delta\phi_n(P) & \Delta\psi_1(P_n) & \cdots & \Delta\psi_m(P_n) \\
\\
\Delta\phi_1(P) & \cdots & 
\cdots & \Delta\phi_n(P) & 0
& \Delta\psi_1(P) & \cdots & \Delta\psi_m(P) \\
\\
\phi_1(Q_1) & \cdots & \cdots & \phi_n(Q_1)
& \psi_1(P) & 0 &  \cdots & \psi_m(Q_1) \\
\\
\vdots & \cdots & \cdots & \vdots & \vdots  & \vdots & \ddots & \vdots \\
\\
\phi_1(Q_m) & \cdots & \cdots & \phi_n(Q_m)
& \psi_m(P) & \psi_1(Q_m) & \cdots & 0 \\
\\

\end{array} \right)
$$

Observe that in this case $K_{N+1}=V(P_{n+1})$
since $\psi_k(P_{n+1})=\phi_{n+1}(Q_k)$ and 
$\Delta\phi_j(P_i)=\Delta\phi_i(P_j)$.

Concerning the determinants, applying Laplace determinantal rule on the last row of $U(P)$ we see that for every $\ell$, $1\leq \ell\leq m$, we get the representation 
\begin{equation} \label{F}
F(P)=\mbox{det}(U(P))=\delta_{N-1}\psi_\ell^2(P)+A(P)\psi_\ell(P)+B(P)
\end{equation}
where 
$$
|\delta_{N-1}|=|\mbox{det}(K_{N-1}(\{P_i\},\{Q_h\}_{h\neq \ell}))|
$$
$$A\in \mbox{span}\{\phi_j,\Delta\phi_j,\psi_k\,;\,1\leq j\leq n\,,\,1\leq k\leq m\,,\,k\neq \ell\}
$$
$$
B\in \mbox{span}\{\phi_i\Delta\phi_j,\psi_k\phi_i,\psi_k\Delta\phi_i,\psi_k\psi_h\,;\,1\leq i,j\leq n\,,\,1\leq k,h\leq m\,,\,k,h\neq \ell\}\;.
$$

Similarly, developing $det(V(P))$ by the $(n+1)$-row we have
\begin{equation} \label{G}
G(P)=\mbox{det}(V(P))= -\mbox{det}(K_{N-1})(\Delta\phi_n(P))^2+C(P)\Delta\phi_n(P)+D(P)
\end{equation}
where $$C\in \mbox{span}\{\Delta\phi_j,\psi_k,\Delta\psi_k\,;\,1\leq j\leq n-1\,,\,1\leq k\leq m\}
$$
$$
D\in \mbox{span}\{\Delta\phi_i\Delta\phi_j,\Delta\phi_i\Delta\psi_h,\psi_k\Delta\phi_i,\psi_k\Delta\psi_h\,;\,1\leq i,j\leq n-1\,,\,1\leq k,h\leq m\}\;.
$$

First, we prove that $G$ is not identically zero in $\Omega$ if $\mbox{det}(K_{N-1})\neq 0$ (the latter a.s. holds by inductive hypothesis). 
Let $P(t)=P_n+t(1,0)$, $t\in\mathbb{R}$, and $r(t)=\|P(t)-P_n\|_2=|t|$. If $G\equiv 0$ then $G(P(t))\equiv 0$ in neighborhood of $t=0$. Then, we would locally have 
\begin{equation} \label{u(t)}
u^2(t)=c(t)u(t)+d(t)\;,\;\;u(t)=\Delta\phi_n(P(t))\;,
\end{equation}
where $c(t)=C(P(t))/\mbox{det}(K_{N-1})$ and $d(t)=D(P(t))/\mbox{det}(K_{N-1})$. Notice that both $c$ and $d$ are analytic in a neighborhood of $t=0$, since $C$ and $D$ are analytic in a neighborhood of $P_n$. By (\ref{u(t)}) and (\ref{lapl}) we get 
\begin{equation} \label{u(t)-2}
u(t)=4\nu t^{2(\nu-1)}\left(\nu\log(|t|)+1\right)\;.
\end{equation}
Clearly $c$ cannot be identically zero there, otherwise $u^2$ would be analytic at $t=0$ and thus would have an algebraic order of infinitesimal as $t\to 0$, whereas by (\ref{u(t)-2}) we have $u^2(t)\sim 16\nu^4t^{4(\nu-1)}\log^2(|t|)$. Hence taking the  Maclaurin expansion of $c$ we get $c(t)\sim c_st^s$ as $t\to 0$ for some $s\geq 0$, the order of the first nonvanishing derivative at $t=0$. Now, $u^2(t)\sim 16\nu^4t^{4(\nu-1)}\log^2(|t|)$, whereas by 
$u^2\equiv cu+d$ we would have $u^2(t)\sim 4\nu^2c_st^{s+2(\nu-1)}\log(|t|)+d_pt^p$, where 
either $d(0)\neq 0$ and $p=0$, or $d(0)=0$ and $p>0$ (the order of the first nonvanishing
derivative at $t=0$). Then   we get a contradiction, since $u^2$ cannot have two distinct limits or orders of infinitesimal at the same point.

Moreover, $G$ is clearly continuous in $\Omega$ and analytic in $\Omega\setminus \{P_1,\dots,P_n\}$, 
since all the functions involved in its definition (\ref{G}) are analytic up to their own center. 
Consequently, if $\mbox{det}(K_{N-1})\neq 0$ by continuity $G$ is not identically zero also in $\Omega\setminus \{P_1,\dots,P_n\}$. 

Then, 
$\mbox{det}(K_{N+1})=\mbox{det}(V(P_{n+1}))=G(P_{n+1})$ is a.s. nonzero, 
since the zero set of a not identically zero real analytic function on an open connected set in $\mathbb{R}^d$ is a null set (cf. \cite{M20} for an elementary proof). More precisely, 
denoting by $Z_G$ the zero set of $G$ in $\Omega$, 
we have that $$Z_G=(Z_G\cap \{P_1,\dots,P_n\})\cup 
(Z_G\cap (\Omega\setminus \{P_1,\dots,P_n\}))\;.
$$
Hence $Z_G$ is a null set if $G\not\equiv 0$, because the first intersection is a finite set, and the second is the zero set of a not identically zero real analytic function. Considering the probability of the corresponding events and recalling that $\mbox{det}(K_{N-1})\neq 0$ (which a.s. holds) implies $G\not\equiv 0$, we can then write 
$$ 
\mbox{prob}\{\mbox{det}(K_{N+1})=0\}=\mbox{prob}\{G(P_{n+1})=0\}
$$
$$
=\mbox{prob}\{G\equiv 0\}
+\mbox{prob}\{G\not\equiv 0\;\&\;P_{n+1}\in Z_G\}
=0+0=0\;,
$$ 
and this branch of the inductive step is completed.

We turn now to the branch of the inductive step where a boundary point is added. In this case we consider the function $F$ in (\ref{F}) restricted to the boundary, that is $F(P(t))$ with $P(t)=\gamma(t)$, $t\in (a,b)$, which for every fixed $\ell\in \{1,\dots,m\}$ has the representation 
$$
F(\gamma(t))=\mbox{det}(U(\gamma(t)))=\delta_{N-1}v^2(t)+A(\gamma(t))v(t)+B(\gamma(t))
$$
where 
\begin{equation} \label{v(t)}
v(t)=\psi_\ell(\gamma(t))=r_\ell^{2\nu}(t)\log(r_\ell(t))\;,\;\;r_\ell(t)=\|\gamma(t)-Q_\ell\|_2
\end{equation}
with $Q_\ell=\gamma(t_\ell)$, $\;t_\ell\in (a,b)$. We claim that if $\delta_{N-1}\neq 0$ (which a.s. holds by inductive hypothesis), $F\circ \gamma$ cannot be identically zero in any of the two connected components of $(a,b)\setminus\{t_1,\dots,t_m\}$ (i.e., the subintervals) having $t_\ell$ as extremum. Otherwise,  
we would have in a left or in a right neighborhood of $t_\ell$ 
\begin{equation} \label{v2(t)}
v^2(t)=\alpha(t)v(t)+\beta(t)\;,
\end{equation}
where $\alpha(t)=A(\gamma(t))/\delta_{N-1}$ and $\beta(t)=B(\gamma(t))/\delta_{N-1}$ are both analytic in a full neighborhood of $t_\ell$. Notice that, since $\gamma'(t_\ell)\neq (0,0)$ (the curve is regular), $r_\ell(t)\sim \|\gamma'(t_\ell)\|_2|t-t_\ell|$ which by (\ref{v(t)}) gives $v(t)\sim \|\gamma'(t_\ell)\|_2^{2\nu}(t-t_\ell)^{2\nu}\l og(|t-t_\ell|)$ and $v^2(t)\sim \|\gamma'(t_\ell)\|_2^{4\nu}(t-t_\ell)^{4\nu}\log^2(|t-t_\ell|)$ as $t\to t_\ell$. Now $\alpha$ cannot be identically zero in any left or right neighborhood, otherwise $v^2\equiv \beta$ there and would have an algebraic order of infinitesimal at $t_\ell$. Hence taking the Taylor  expansion of $\alpha$
we get $\alpha(t)\sim \alpha_s(t-t_\ell)^s$ as $t\to t_\ell$ for some $s\geq 0$, the order of the first nonvanishing
derivative at $t=t_\ell$. 
On the other hand, by 
$v^2\equiv \alpha v+\beta$ locally, we would have $v^2(t)\sim \|\gamma'(t_\ell)\|_2^{2\nu}\alpha_s(t-t_\ell)^{s+2\nu}\log(|t-t_\ell|)+\beta_p(t-t_\ell)^p$, where either $\beta(t_\ell)\neq 0$ and $p=0$, or $\beta(t_\ell)=0$ and $p>0$ (the order of the first nonvanishing
derivative at $t=t_\ell$). Again we get a contradiction, since $v^2$ cannot have two distinct limits or orders of infinitesimal at the same point.

The result is that $F\circ \gamma$ is a.s. not identically zero in any 
connected component of $(a,b)\setminus \{t_1,\dots,t_m\}$. 
Then, 
$\mbox{det}(K_{N+1})=\mbox{det}(U(Q_{m+1}))=F(\gamma(t_{m+1}))$ is a.s. nonzero. In fact, observe that $F\circ \gamma$ is analytic in $(a,b)\setminus \{t_1,\dots,t_m\}$, since $F$ is analytic in $\mathbb{R}^2\setminus (\{Q_1,\dots,Q_m\}\cup \{P_1,\dots,P_n\})$. Moreover, 
denoting by $Z_{F\circ \gamma}$ the zero set of $F\circ \gamma$ in $(a,b)$, 
we have that $$Z_{F\circ \gamma}=(Z_{F\circ \gamma}\cap \{t_1,\dots,t_m\})\cup 
(Z_{F\circ \gamma}\cap ((a,b)\setminus \{t_1,\dots,t_m\}))\;.
$$
Hence $Z_{F\circ \gamma}$ is a null set if ${F\circ \gamma}\not\equiv 0$, because the first intersection is a finite set, and the second is the componentwise finite union of the zero sets of a not identically zero real analytic function on each connected component. Considering the probability of the corresponding events and recalling that $\mbox{det}(K_{N-1})\neq 0$ (which a.s. holds) implies $F\circ \gamma\not\equiv 0$, we can then write 
$$ 
\mbox{prob}\{\mbox{det}(K_{N+1})=0\}=\mbox{prob}\{F(Q_{m+1})=0\}
$$
$$
=\mbox{prob}\{F\circ \gamma\equiv 0\}
+\mbox{prob}\{F\circ \gamma\not\equiv 0\;\&\;t_{m+1}\in Z_{F\circ \gamma}\}
=0+0=0\;,
$$ 
and also the boundary branch of the inductive step is completed.
\hspace{0.2cm} $\square$

\subsection{Remarks on possible extensions}
The result of Theorem 1 is a first step towards a theory of Kansa  collocation unisolvence, and could be extended in several directions within the random framework. The first extension comes immediately from the fact that a null set has also measure zero for any continuous measure with density (that is, absolutely continuous with respect to the Lebesgue measure). We can state indeed the following
\begin{theorem}
The assertion of Theorem 1 holds true if the points $\{P_i\}$ are independent identically distributed with respect any continuous probability measure 
with density on $\Omega$, say $\sigma\in L^1 _+(\Omega)$, and the abscissas $\{t_h\}$
are independent identically distributed with respect any continuous probability measure 
with density on $(a,b) $, say $w\in L^1 _+(a,b)$.
\end{theorem}

This extension could be interesting whenever it is known that the solution has steep gradients or other regions where it is useful to increase the discretization density. Concerning the implementation of random sampling with respect to continuous probability densities, we recall the well-known ``acceptance-rejection method'', cf. e.g. \cite{CCTC22,Flury90,NO16} with the references therein.

More difficult but worth of further investigations are:
\begin{itemize}

\item extension to $\Omega\subset \mathbb{R}^d$, $d\geq 3$;

\item extension to other analitic RBF up to the center, e.g. Radial Powers; 

\item extension to piecewise analytic boundaries;

\item extension to other differential operators and/or boundary conditions.

\end{itemize}
The latter in particular could be challenging, since the operators involved in the equation and in the boundary conditions may not be radial.

\vskip0.5cm 
\noindent
{\bf Acknowledgements.} 

Work partially supported by the DOR funds of the University of Padova, and by the INdAM-GNCS 2024 Projects “Kernel and polynomial methods for approximation and integration: theory and application software''. 

This research has been accomplished within the RITA ``Research ITalian network on Approximation" and the SIMAI Activity Group ANA\&A, and the UMI Group TAA ``Approximation Theory and Applications".

\end{document}